\documentclass [12pt,a4paper,reqno]{amsart}
\usepackage{amsmath,amsthm}
\usepackage{amscd}
\usepackage{amsfonts}
\input amssymb.sty
\DeclareMathOperator{\Hom}{Hom}

\newcommand{\ef}{\end{equation}}
\chardef\bslash=`\\ 

\newtheorem{thm}{Theorem}
\newtheorem*{thm*}{Theorem}
\newtheorem*{conjecture*}{Conjecture}

\newtheorem{cor}{Corollary}

\newtheorem{lem}{Lemma}

\theoremstyle{definition}
\newtheorem{defn}{Definition}
\newtheorem*{remark*}{Remarks}
\newtheorem*{examples*}{Examples}
\newtheorem*{defn*}{Definition}
\newtheorem*{cor*}{Corollary}

\newtheorem{rem}{Remark}

\newcommand{\thmref}[1]{Theorem~\ref{#1}}
\newcommand{\secref}[1]{Section~\ref{#1}}

\newcommand{\lemref}[1]{Lemma~\ref{#1}}
\newcommand{\corref}[1]{Corollary~\ref{#1}}

\numberwithin{equation}{section}

\newcommand{\Th}{\Theta}
\newcommand{\C}{\mathcal{C}}
\newcommand{\D}{\mathcal{D}}
\newcommand{\E}{\mathcal{E}}
\newcommand{\V}{\mathcal{V}}
\newcommand{\Q}{\mathbf{Q}}
\newcommand{\st}{\sigma}

\renewcommand{\sectionmark}[1]{}

\date{}
\begin{document}

\title{On automorphisms of categories of universal algebras}

\author[Boris Plotkin]{B. Plotkin\\
 Institute of Mathematics \\
 Hebrew University, 91803 Jerusalem, Israel\\}

\author[Grigori Zhitomirski]{G. Zhitomirski\\
 Department of Mathematics \\
 Bar-Ilan University, 52900 Ramat Gan, Israel}
\thanks {This research of the second author was partially supported by THE ISRAEL SCIENCE FOUNDATION
 founded by The Israel Academy of Sciences and Humanities - Center of
 Excellence Program.}

\maketitle

\begin{abstract}
Given a variety $\V$ of universal algebras. A new approach is suggested to characterize algebraically
automorphisms of the category of free $\V$-algebras. It gives in many cases an answer to the problem set by the
first of authors, if automorphisms of such a category are inner or not. This question is important for universal
algebraic geometry \cite{Free,AlgGeom}. Most of results will actually be proved to hold for arbitrary categories
with a represented forgetful functor.
\end{abstract}


\baselineskip 20pt
\bigskip

\centerline{INTRODUCTION}\label{Intro}

It is a current opinion that the notions of an isomorphism and an automorphism of categories are not important.
As far as the authors know there are no researches devoted to describing automorphisms of categories although
this theme is very popular for the most of other algebraic structures.  The first of authors set the problem to
describe automorphisms of a category of free algebras of some given variety of universal algebras. It turns out
that this problem is quite important for universal algebraic geometry \cite{Free,AlgGeom}. The most important
case is, when all automorphisms of a category in question are inner or close to inner in a sense.

Recall that an automorphism $\Phi$ of a category $\C$ is called inner if it is isomorphic to the identity
functor $Id \sp {\C}$   in the category of all endofunctors of $\C$. It means that for every object $A$ of the
given category there exists an isomorphism $\sigma \sb A : A\to \Phi (A)$ such that for every morphism $\mu:
A\to B$ we have $\Phi (\mu )= \sigma \sb B \circ \mu \circ \sigma \sb A \sp {-1}$. This fact explains the term
"inner". Thus if an automorphism $\Phi$ is inner the object $\Phi (A)$ is isomorphic to $A$ for every
$\C$-object $A$.

Let $\V$ be a variety of universal algebras. Consider the category $\Th (\V)$ whose objects are all algebras
from $\V$ and whose morphisms are all homomorphisms of them. Fix an infinite set $X\sb 0$. Let $\Th \sp 0 (\V )$
be the full subcategory of $\Th (\V )$ defined by all free algebras from $\V$ over finite subsets of the set
$X\sb 0$. The group of automorphisms of the category $\Th \sp 0 (\V )$ is the subject of inquiry.

It is known, for example, that every automorphism of the category  $\Th \sp 0 (\V )$  is inner if $\V $ is the
variety of all groups \cite{Free}. But it is not so, for example, if $\V $ is the variety of all semigroups
\cite{Free} or the variety of all Lie algebras \cite{LieAlg}. In the last cases automorphisms are similar to
inner but not inner exactly. What does it mean "similar to inner"? It means different things for different
varieties. For instance, if $\V $ is the variety of all semigroups the maps $f\sb A : A\to \Phi (A)$ mentioned
above are not necessarily isomorphisms but can be anti-isomorphisms too.  If $\V $ is the variety of all modules
over a ring $\bf K $, these maps are so called semi-automorphisms, more exactly they are pairs of maps
$(f,\sigma )$, where $f$ is an automorphism of the ring $\bf K $, $\; \sigma :M\to N$ is an additive bijection
satisfying the condition: $\sigma (ax)=f(a)\sigma (x)$ for all $a\in K$ and $x\in M$.

There exists a usual approach to the question if all automorphisms of the category $\Th \sp 0 (\V )$  are inner
or semi-inner. It demands describing the group $AUT\;END(F)$ of all automorphisms of the semigroup of all
endomorphisms of the free algebras $F$ of the given variety $\V$, and then applying the so called Reduction
theorem (\thmref{Plotkin}) due to B. Plotkin \footnote {It should be mention that this theorem was first proved
by Berzins \cite{Berz} for the variety of commutative associative algebras over an infinite field.}. This
theorem gives an opportunity to make a conclusion about automorphisms of the category $\Th \sp 0 (\V )$ (under
some conditions) if all automorphisms of $ END(F)$ are described, where finitely generated free algebra $F$
generates the variety $\V$. The original proof of this theorem is based on several special algebraic notions and
constructions.

But it turns out that going over to a general case of categories supplied with a represented forgetful functor
we can solve the mentioned problem as a whole and obtain very easily this theorem and even some more
applications. And what is more a new approach can be suggested to answer the mentioned question. To describe the
group $AUT\;END(F)$ for a given free algebra $F$ is a not trivial problem (see for example \cite{Form,InvSem}).
Instead of this it is sufficient to find out how an automorphism of the given category acts on the set of all
homomorphisms from a monogenic free algebra into a finitely generated free algebra. And the last problem is
reduced to a purely algebraic problem to study some derivative (polynomial) operations in a free algebra, in
most cases in a two-generated free algebra. This is the purpose of the present paper.

First of all a notion of {\it potential-inner} automorphisms is introduced (\secref{sect_inner}). Let $\C $ be a
subcategory of a category $\D$.  An automorphism $\Phi$ of the category $\C$ is called $\D-inner$ if for every
object $A$ of the category $\C$ there exists an $\D-$isomorphism $f\sb A : A\to \Phi (A)$ such that $\Phi (\mu)=
f\sb B \circ \mu \circ f\sb A \sp {-1}$ holds for every $\C$-morphism $\mu: A\to B$. An automorphism $\Phi $ of
a given category $\C$ is called potential-inner if it is $\D$-inner for some extension $\D$ of $\C$. This notion
gives an opportunity to consider the problem from a new point of view. The necessary and sufficient condition
found for an automorphism to be potential-inner (\thmref{potent-inner}) is satisfied in all important cases.

Thus the main problem can be formulated now in the following way: 1) what extension $\D$ of the given category
$\C$ we have to construct in order to make all $\C$-automorphisms to be $\D$-inner and 2) when potential-inner
automorphisms are in fact inner. These problems are solved in \secref{sect_Univ} for categories of free algebras
using the derivative operation language (\lemref{InnerAndCentral}, \thmref{main2}). Roughly speaking, all
potential-inner automorphisms are produced by isomorphisms of algebras onto derived algebras, a potential-inner
automorphism is in fact inner, if there exist so called central isomorphisms of given algebras onto derived
algebras (algebras on the same underlying sets with respect to derived operations).

Secondly, notions are introduced of left and right {\it indicators} in a category (\secref{sect_Reduct}). Let
$\Q: \C \to Set$ be a forgetful functor. An object $A\sp 0$ of the category $\C$ is called a {\it right
indicator} if for every two objects $A$ and $B$ and for every bijection $s: \Q (A )\to \Q (B )$ the following
condition is satisfied:

if for every morphism $\nu: B \to A\sp 0  $  there exists a morphism $\mu: A\to A \sp 0 $ such that $ \Q (\mu )=
\Q (\nu )\circ s $, then there exists an isomorphism $\gamma : A \to B$ such that $\Q (\gamma )=s$.

Dually the notion of a left indicator is defined.

For example, every free algebra generating a variety $\V$ is a right indicator in the category $\Th (\V)$ and
every free algebra $H$ in $\V$ having not less free generators than arities of all operations is a left
indicator. Using these notions, a generalized reduction theorem is proved (\thmref{main}, \thmref{general}) and
it is shown in \secref{sect_Univ} how the original Reduction Theorem can be obtained from it (\thmref{Plotkin}).

Further in this section, the method mentioned above is described and applied to categories of semigroups
(\thmref{semigroup}) and inverse semigroups (\thmref{inverse}) only to show how it works. In our next paper
\cite{Pl_Zh}, we apply this method and characterize automorphisms of categories  $\Th \sp 0 (\V) $ in the case
$\V $ is the variety of all associative $K-$algebras, where $K$ is a infinite field, and in the case $\V $ is
the variety of all group representations in unital $R-$ modules, where $R$ is an associative commutative ring
with unit.

In the last section (\secref{applications}) some applications are given for categories of sets and semigroups of
transformations  (\thmref{set}, \thmref{Schreier}, \thmref{Shutov}).

For the some notions and results of Category Theory and Universal Algebra we refer a reader to \cite{Categ,
Gratz}. Some part of results presented in this paper was published by the second author in \cite{Zh}.

{\bf Acknowledgements.} The authors are happy to thank E.Plotkin and G.Mashevitzky for the stimulating
discussions of the results.

\section{Inner and potential-inner automorphisms}\label{sect_inner}
We consider such categories $\C$ that are represented in the category $Set$ (the category of all sets and maps),
that is, there exists a faithful functor  $\Q: \C \to Set$. Such a functor is called a forgetful functor. If
$\C$ is a category of universal algebras, then the forgetful functor is usually the natural forgetful functor,
which assigns to every algebra $A$  the underlying set $\vert A \vert$ and to every homomorphism itself as a
mapping, but not only this case.

We assume that a forgetful functor for every category we consider is fixed, and we say that a category $\C$ is a
subcategory of a category $\D$ having in the mind that the forgetful functor for the category $\D $ is an
extension of the forgetful functor for the category $\C$.

If $\C $ is a subcategory of a category $\D$ then $Id\sp {\D } \sb {\C}: \C \to \D $ denotes the natural
embedding functor, that is the restriction of identity functor $Id \sp {\D}$ to $\C $.

\begin{defn} Let $\C $ be a subcategory of a category $\D$.  A functor $\Phi : \C \to \D $ is said to be {\it inner}
if it is isomorphic in the sense of category of functors to the functor $Id\sp {\D } \sb {\C}$.
\end{defn}
It is clear that every inner functor is faithful.
\begin{defn} Let $\C $ be a subcategory of a category $\D$.  An automorphism $\Phi$ of the category $\C$ is
called $\D-inner$ if the functor $Id\sp {\D } \sb {\C}\circ \Phi :\C \to \D $ is inner. That is for every object
$A$ of the category  $\C$ there exists an $\D-$isomorphism $f\sb A : A\to \Phi (A)$ such that for every $\C
-$morphism $\mu: A\to B$ we have $\Phi (\mu)= f\sb B \circ \mu \circ f\sb A \sp {-1}$. That is, the following
diagram is commutative:
$$
\CD
A @>f\sb A>>    \Phi (A)\\
@V\mu VV          @VV \Phi (\mu ) V\\
B @>f\sb B  >>   \Phi (B).
\endCD
$$
A $\C$-inner automorphism of a category $\C$ is called an {\it inner} automorphism. A category is called {\it
perfect} if all its automorphisms are inner.
\end{defn}

In other words, an automorphism of the category $\C $ is $\D -$ inner if it can be extended to an inner
automorphism of the category $\D$.

If an automorphism $\Phi $ of $\C$ is $\D -$inner , a family $(f\sb A : A\to \Phi (A) \, \vert A\in Ob \, \C)$
exists with the condition above. But this family is not only one possible. For example, the identity
automorphism is of course inner, because it is determined by the family of identities, but it may be there
exists another family of $\D -$morphisms which determines it.

\begin{defn}\label{central} A function $c$ which assigns a permutation $c\sb A :\Q (A)\to \Q (A) $ to every
object $A$ of $\C$ is called {\it a central function} if it determines the identity automorphism of $\C$, that
is, $c\sb B \circ \Q (\mu )\circ c\sb A \sp {-1} =\Q (\mu )$ for every $\C -$morphism $\mu: A\to B$.
\end{defn}

It is obvious that two families $(f\sb A :A\to \Phi (A) \vert A\in Ob \, \C)$ and $(g\sb A :A\to \Phi (A) \vert
A\in Ob \, \C)$ of $\D -$isomorphisms determine the same automorphism $\Phi $ of the category $\C$ if and only
if the corresponding maps are equal up to a central function: $\Q (f\sb A)=\Q (g\sb A ) \circ c\sb A$ for all
objects $A$.

\begin{defn} Let $\C$ be a category with a forgetful functor $\Q: \C \to Set$. An automorphism $\Phi$ of the
category $\C$ is said to be  {\it potential-inner} if it is $\D$-inner for some category $\D$ such that $\C$ is
a subcategory of $\D$ with the same objects.
\end{defn}

We illustrate the last definition with the following examples.

 {\bf  Examples.}

1. Let $\C $ be the category of semigroups and homomorphisms, in this case $\D $ can be the category of
semigroups and homomorphisms and anti-homomorphisms, the functor $\Q$ is the natural forgetful functor. Thus in
the definition above $f\sb A $ can be either an isomorphism or an anti-isomorphism.

2. Let $\C$ be a category of unitary modules over a ring $\bf K$. The functor $\Q$ is the natural forgetful
functor. The category $\D$ has the same objects but its morphisms are additives maps $\sigma $ of modules $M$ to
$N$ satisfying the condition: $\sigma (ax)=f(a)\sigma (x)$ for all $a\in K$ and $x\in M$, where $f$ is an
automorphisms of the ring $\bf K$. The category $\C$ can be identified with the subcategory of $\D$ for which
$f=1\sb {\bf K}$.

3. (See \cite{SevLect}). Let $\C$ be an arbitrary category of algebras and $G$ be a fixed non-trivial algebra in
$\C$. Denote by $\C \sp G $ the category which objects are $\C$-monomorphisms with the same domain $G$  and
morphisms of which are commutative diagrams:
$$
\CD
G @>h>>    H\\
@V1\sb G VV     @VV \mu V\\
G @>h'>>    H',
\endCD
$$
where $h,h'$ are monomorphisms and  $\mu $ is a morphism in $\C$. The objects of this category can be considered
as $G$-algebras $H$, that is, the algebras with fixed algebra of constants $G$. In this way, to elements $g\in
G$ correspond constants, i.e., nullary operations in $H$.

The category $\C \sp G $ is a subcategory of the category $\C (G)$ that has the same objects and the morphisms
of which are the commutative diagrams:
$$
\CD
G @>h>>    H\\
@V\sigma VV     @VV \mu V\\
G @>h'>>    H',
\endCD
$$
where $\sigma $ is now an automorphism of $G$.

Define a functor $\Q$ in the same way as in the previous example: $\Q (h:G\to H)= \vert H\vert$ and $\Q$ assigns
to every commutative diagram above the map $ \mu : \vert H\vert \to \vert H'\vert$. Let the category $\D$ in
Definition~2  be $\C (G)$, then $\D$-inner automorphisms of the category $\C \sp G $ are  {\it semi-inner }
automorphisms of this category introduced in \cite{SevLect}.

\vspace{12pt}
If all automorphisms of a category $\C$ are potential inner the problem is to find the smallest
its extension $\D$ such that all automorphisms of $\C$ are $\D -$inner. We start with several very simple facts.
\begin{lem}\label{subcategory}
Let $\C $ be a subcategory of a category $\D$ and $\Phi :\C \to \D $ be a functor. Let $\E $ be a subcategory of
$\C$ such the restriction of $\Phi $ to $\E $ is a $\D$-inner automorphism of $\E $. Then $\Phi $ is a
composition of two functors $\Phi = \Psi \circ \Gamma $, where $\Gamma :\C \to \D $ is an identity on $\E $
(preserves all objects and all morphisms of $\E$), and the functor $\Psi$ is an inner automorphism of the
category $\D$.

\end{lem}
\begin{proof}
Let the family $(s\sb A : A\to \Phi (A)\;\vert A\in Ob\E)$ be an isomorphism of functors $Id\sp {\D } \sb {\E}
\to \Phi\sb {\vert \E }$, that is, $\Phi (\nu )=s\sb B \circ \Phi (\nu )\circ s\sb A \sp {-1}$ for every
$\E$-morphism $\nu :A\to B$.  We construct an inner automorphism  $\Psi $ of $\D $ in the following way. For
every $\D $-object $X$ , we set $\Psi (X)=X $ if $X$ is not an object of $\E $ , and $\Psi (A) =\Phi (A) $ for
every object $A$ of $\E $. Further, we define $\D$-isomorphism $u\sb X :X\to \Psi (X)$ in the following way:
$u\sb X =1\sb X$ if $X $ is not an object of  $\E $, and $u\sb A = s\sb A $, where $A\in Ob\, \E $. For every
$\D$-morphism $f:X\to Y$, let $\Psi (f) = u\sb Y \circ f \circ u\sb X \sp {-1}$. According to the given
construction, $\Psi $ is an inner automorphism of $\D $.

It is clear that $\Phi = \Psi \circ \Psi \sp {-1} \circ \Phi $. Let $\Gamma =\Psi \sp {-1} \circ \Phi $.
According to this definition, we have that $\Gamma (A)=(\Psi \sp {-1} \circ \Phi )(A)=A $ for all $A\in Ob\,\E$
and $\Gamma (\nu )=(\Psi \sp {-1} \circ \Phi )(\nu)=s\sb B \sp {-1}\circ \Phi (\nu )\circ s\sb A =\nu $ for all
$\E$-morphisms $\nu :A\to B$.
\end{proof}

Assuming in the above lemma a category $\E $ to be discrete, that is, it contains no arrows besides identities,
we obtain as a consequence the following result (see also \cite{LieAlg}).

\begin{lem}\label{objects}
Let  $\Phi $ be an automorphisms of a category $\C$. Suppose further that for some class $\bf E$ of $\C
-$objects,  $\Phi (A)$ is isomorphic to $A$ for every $A\in \bf E$. Then $\Phi$ is a composition of two $\C
-$automorphisms $\Phi = \Psi \circ \Gamma $, where $\Gamma $ leaves fixed all objects from $\bf E$ and $\Psi$ is
an inner automorphism.

\end{lem}
\begin{proof}
If the class $\bf E$ is closed under $\Phi$ and $\Phi \sp {-1}$, apply the previous lemma to the subcategory $\E
$ whose class of objects is $\bf E$ and whose morphisms are identity morphisms only (discrete subcategory) and
obtain the required result where $\D =\C$. If $\bf E$ is not closed we can consider its $\Phi -$ and $\Phi \sp
{-1}-$closure that has clearly the same property.
\end{proof}

The next result will give us a necessary condition for an automorphism to be potential inner.
\begin{lem}\label{free objects}
Let $F$ be a free object in a category $\C$ over a set $X$. If an automorphism $\Phi $ of $\C$ is potential
inner then $\Phi (F)$ is also a free object over $X$ and hence it is isomorphic to $F$.
\end{lem}
\begin{proof}
Under hypothesis, there is a map $m:X\to \Q (F)$ such that for every $\C -$object $A$ and for every map $f:X\to
\Q (A)$, there exists an unique morphism $\bar{f}:F\to A$ such that $f=\Q (\bar{f})\circ m$. If the given
automorphism $\Phi $ of $\C$ is potential inner, that is, it is $\D -$inner for some extension $\D$ of the
category $\C$, then there exist a family of $\D -$isomorphisms $s\sb A :A \to \Phi (A)$ such that $\Phi (\nu
)=s\sb B \circ \nu \circ s\sb A \sp {-1}$ for every $\C -$morphism $\nu :A\to B$. Set $\tilde{m}=\Q (s\sb F
)\circ m$. It is a map from $X$ to $\Q (\Phi (F))$. Let $A$ be an object and $f: X\to \Q (A)$ a map. Consider
the object $B=\Phi \sp {-1} (A)$ and the map $g=\Q (s\sb B \sp {-1})\circ f$ from $X$ to $B$. Then we have an
unique morphism $\bar{g}:F\to B$ such that $g=\Q (\bar{g})\circ m$. Hence we have a morphism $\nu =\Phi
(\bar{g})$ from $\Phi (F) $ to $A$ with the following condition:
$$\Q (\nu )\circ \tilde{m}=\Q ( s\sb B \circ \bar{g} \circ s\sb F \sp {-1})\circ \Q
(s\sb F )\circ m =\Q(s\sb B \circ \bar{g} \circ s\sb F \sp {-1}\circ s\sb F)\circ m =\Q(s\sb B )\circ \Q
(\bar{g})\circ m =$$
$$=\Q(s\sb B )\circ g=\Q(s\sb B )\circ \Q (s\sb B \sp {-1})\circ f =f.$$
The uniqueness if the morphism $\nu$ with the condition $\Q (\nu )\circ \tilde{m}=f$ is clear. Thus $\Phi (F)$
is also a free object over $X$ with respect to the map $\tilde{m} :X\to \Q (\Phi (F))$.
\end{proof}
We consider categories $\C$ such that the forgetful functor $\Q :\C \to Set$ is represented by a pair $(A\sb 0 ,
x\sb 0 )$ where $A\sb 0 $ is an object of $\C$ and $x\sb 0 \in \Q (A\sb 0 )$. It means that for every element
$a\in \Q (A)$ for some object $A$ there exists an unique $\C$-morphism $\alpha \sb a \sp A : A\sb 0 \to A$ such
that

\begin{equation}\label{alphas}
\Q (\alpha \sb a \sp A )(x\sb 0)=a,
\end{equation}
in other words $A\sb 0 $ is a free object over set $\{x\sb 0 \}$. In the case the object $A$ is known, we turn
down the letter $"A"$ in the designation $\alpha \sb a \sp A $.

If an automorphism $\Phi $ of $\C $ is potential inner, then, according to \lemref{free objects}, $\Phi (A\sb
0)$ is isomorphic to $A\sb 0 $. If $\Phi (A\sb 0)$ is isomorphic to $A\sb 0 $, then (according to
\lemref{objects}) $\Phi $ is a composition of two automorphisms $\Phi = \Psi \circ \Gamma $, where $\Gamma $
preserves $A\sb 0 $ , and the functor $\Psi$ is a $\D$-inner automorphism of the category $\C$. Therefore we can
restrict ourself to considering automorphisms that preserve the object $A\sb 0$.

Let $\Phi $ be such an automorphism. Consider an arbitrary object $A$. It was mentioned above that there is a
bijection $a\mapsto \alpha \sb a \sp A$ between sets $\Q (A)$ and $\Hom (A\sb 0 ,A)$ defined by \ref{alphas}.
Define a map  $s \sb A \sp {\Phi }:\Q (A)\to \Q (\Phi(A))$ setting for every $a\in \Q (A) $:
\begin{equation}\label{E:sDefinition}
 s \sb A \sp {\Phi }(a) =\bar{a}\Leftrightarrow \Phi (\alpha \sb a \sp A)=\alpha \sb {\bar {a}}\sp {\Phi (A)}
\end{equation}

or
\begin{equation}\label{E:sDefDirect}
 s \sb A \sp {\Phi }(a) =\Q (\Phi (\alpha \sb a \sp A))(x\sb 0).
\end{equation}

It is very simple to verify that the function $s\sp {\Phi }: A\mapsto s \sb A \sp {\Phi }$ has the following
properties:
\begin{equation}
s \sb A \sp {Id}=1\sb A,
\end{equation}
\begin{equation}
s \sp {\Phi \sp {-1}}=(s \sp {\Phi })\sp {-1},
\end{equation}
and for two automorphisms $\Phi $ and $\Psi$
\begin{equation}
s \sp {\Psi \circ \Phi }=s \sp {\Psi }\circ s\sp {\Phi }.
\end{equation}

If an automorphism $\Phi $ is fixed, we do not write the superscript $\Phi$. For every morphism $\nu : A \to B$
and for every $a\in \Q (A)$, we have $\nu \circ \alpha \sb a \sp A= \alpha \sb {\Q (\nu )(a)}\sp B$. Applying to
this equation the automorphism $\Phi $ we obtain: $\; \Phi (\nu )\circ \alpha \sb {s \sb A (a)}\sp {\Phi (A)}
=\alpha \sb {s\sb B (\Q (\nu )(a))}\sp {\Phi (B)}\;$. Now apply the functor $\Q $ and take the common value of
corresponding maps in the point $x\sb 0$. We have $\Q (\Phi (\nu ))\circ s \sb A (a)=s\sb B (\Q (\nu )(a))$.
Since $a$ is an arbitrary element of $\Q (A)$ we obtain finally:
 $\; \Q (\Phi (\nu ))\circ s \sb A = s \sb B \circ \Q (\nu )\;$ and hence
 \begin{equation}\label{E:general}
\Q (\Phi (\nu ))=s \sb B \circ \Q (\nu )\circ s \sb A \sp {-1}.
\end{equation}

Notice that the map $s\sb {A\sb 0}$ satisfies some special condition: $s\sb {A\sb 0} (x\sb 0 )=x\sb 0$. We
formulate the obtained result in the following lemma.
\begin{lem}\label{mainlem}
Let $\C$ be a category with a forgetful functor $\Q: \C \to Set$  that is represented by a pair $(A\sb 0 ,x\sb
0)$. If $\Phi $ is an automorphism of the category $\C $ that leaves fixed the object $A\sb 0 $ then there
exists a family of bijections $(s\sb A :\Q (A)\to \Q(\Phi (A)),\; \vert A\in Ob\C )$, such that for every
$\C$-morphism $\nu :A\to B$ we have:

$$\Q (\Phi (\nu ))=s \sb B \circ \Q (\nu )\circ s \sb A \sp {-1}$$
and $s\sb {A\sb 0} (x\sb 0 )=x\sb 0$.

Further, if $\D $ is an extension of $\C$  and the functor $\Q $ can be extended to a functor from $\D$ to
$Set$, such that for every $\nu: A\to B$
$$
\Phi (\nu )=\sigma \sb B \circ \nu \circ \sigma \sb A \sp {-1},
$$
where $\sigma \sb A :A\to \Phi (A)$ are $\D$-isomorphisms  and $\Q (\sigma \sb {A\sb 0})(x\sb 0 )=x\sb 0$, then
$$\Q (\sigma \sb A )=s\sb A $$ for every $\C$-object $A$.
\end{lem}
\begin{proof}
Let the pair $(A\sb 0 , x\sb 0)$ represent the functor $\Q $.  The first statement of the lemma is proved above.
The second one follows from the hypotheses immediately:  $s\sb A (a)=\Q (\Phi (\alpha \sb a ))(x\sb 0 )=\Q
(\sigma \sb A) \circ \Q (\alpha \sb a )(x\sb 0 )=\Q (\sigma \sb A )(a)$ for every $a\in \Q (A)$, that is, $s\sb
A =\Q (\sigma \sb A )$.

\end{proof}

\begin{thm}\label{potent-inner} Let $\C$ be a category with a forgetful functor $\Q: \C \to Set$  that is represented
by an object $A_0$. An automorphism $\Phi$ of the category $\C$ is potential-inner if and only if the $A\sb 0 $
and $\Phi (A\sb 0)$ are isomorphic.
\end{thm}
\begin{proof}
The necessity of that condition follows from \lemref{free objects}. We prove that it is sufficient. Under
hypothesis, the category $\C$ is isomorphic to a subcategory of the category $Set$, and we can assume that $\C$
is a subcategory of $Set$ with the trivial forgetful functor $Id\sp {Set} \sb {\E}$.  Now we can assume that the
automorphism $\Phi$ preserves $A\sb 0 $. According to \lemref{mainlem} there exists a family of bijections
$(s\sb A : A\to \Phi (A),\; \vert A\in Ob\C )$, such that for every $\C$-morphism $\nu :A\to B$ we have:
$$ \Phi (\nu )=s \sb B \circ \nu \circ s \sb A \sp {-1}.$$
Adding these bijections $s\sb A$ and their inverses $s \sb A \sp {-1}$ to the category $\C $ we obtain a new
category of sets $\D$ containing $\C$ as a subcategory with the same objects. Under definition the automorphism
$\Phi$ is $\D$-inner.

\end{proof}

We assume now that automorphisms of the category $\C$ we consider leave fixed the object $A\sb 0$. According to
\lemref{mainlem}, describing of such automorphisms is reduced to the problem of finding out what maps  $s\sb
{A}$ are. The following simple facts are very useful. First we obtain a corollary from \lemref{subcategory}.

\begin{cor}\label{help}
Let $A$ be an object. Suppose that $\Phi (A) =A$ and for some subset $X \subseteq \Q (A) $ there exists an
automorphism $\sigma $ of $A$ such that $\Phi (\alpha \sb x \sp A) =\alpha \sb {\sigma (x)}\sp A$ for every
$x\in X$ (in other words, $s\sb A \sp {\Phi }(x)=\sigma (x)$). Then $\Phi$ is a composition of two automorphisms
$\Phi =\Psi \circ \Gamma $, where $\Psi $ is an inner automorphism and $\Gamma (\alpha \sb x \sp A)= \alpha \sb
x \sp A$ for all $x\in X$, that is $s\sb A \sp {\Gamma }(x)=x$.
\end{cor}
\begin{proof}
Consider a subcategory $\E$ of $\C $ whose objects are $A\sb 0 $ and $A$ only and whose morphisms are identities
on these two objects and all morphisms $\alpha \sb x \sp A$ for $ x\in X$ only. Under hypotheses, for every
$x\in X$ we have: $\Phi (\alpha \sb x \sp A) =\sigma \circ \alpha \sb x \sp A \circ 1\sb {A\sb 0}$. It means
that $\Phi $ induces an inner automorphism of $\E$. According to \lemref{subcategory}, $\Phi$ is a composition
of two automorphisms $\Phi =\Psi \circ \Gamma $, where $\Psi $ is an inner automorphism and $\Gamma $ is an
identity on $\E$, that is, $\Gamma (\alpha \sb x \sp A)= \alpha \sb x \sp A$ for all $x\in X$.
\end{proof}

Now we apply this general fact to a special situation.

\begin{cor}\label{generators}
Suppose that $A=F(X)$ is a free object in $\C $ over a set $X$ and $m:X\to \Q (A)$ the corresponding map.
Suppose that $\Phi (A) =A$. Then $\Phi$ is a composition of two automorphisms  $\Phi =\Psi \circ \Gamma $, where
$\Psi $ is an inner automorphism and $\Gamma (\alpha \sb {m(x)}\sp A )= \alpha \sb {m(x)}\sp A $ for all $x\in
X$, that is $s\sb A \sp {\Gamma }(m(x))=m(x)$.
\end{cor}
\begin{proof}
Denote for simplicity $s=s\sb A \sp {\Phi }$.  Consider two maps $s\circ m:X\to \Q (A)$ and $s\sp {-1} \circ m
:X\to \Q (A)$. They give us two morphisms $\sigma :A\to A $ and $\tau :A\to A$ such that $\Q (\sigma )\circ m
=s\circ m $ and $\Q (\tau ) \circ m =s\sp {-1} \circ m $. These both conditions can be written in the following
way:
\begin{equation}\label{*}
\sigma \circ \alpha \sb {m(x)}=  \alpha \sb {s(m(x))} ,\; \tau \circ \alpha \sb {m(x)}=  \alpha \sb {s\sp {-1}
(m(x))}
\end{equation}
for all $x\in X$. Applying $\Phi \sp {-1}$ to the first one, we obtain $\Phi \sp {-1} (\sigma )\circ \alpha \sb
{s\sp {-1} (m(x))}=\alpha \sb {m(x)}$ and hence $\Q (\Phi \sp {-1} (\sigma ))(s\sp {-1} (m(x))=m(x)$ for all
$x\in X$. It means that
$$\Phi \sp {-1} (\sigma )\circ \tau =1\sb A
$$
and hence
$$ \sigma \circ \Phi (\tau )=1\sb A .
$$
Applying $\Phi $ to the second condition in \ref{*}, we obtain in the same way that
$$ \Phi (\tau )\circ \sigma  =1\sb A .
$$
These equalities give the fact that $\sigma $ and $\tau $ are inverse, that is they are automorphisms of $A$.
Applying \corref{help} we obtain the required result.
\end{proof}

The last result gives us opportunity to restrict our consideration to the case of automorphisms $\Phi $ such
that $s\sb A \sp {\Phi }(m(x))=m(x)$ for a given free object $A$. The next result shows how it simplifies the
situation. Denote by $\theta \sb f$ the unique endomorphism of the object $A$ such that $\Q (\theta \sb f)
(m(x))=f(x)$ for all $x\in X$, where $f:X\to \Q (A)$ is a given map.

\begin{cor}\label{thetas} $\Phi (\theta \sb f) =\theta \sb {s\sb A \sp {\Phi}\circ f } $.
\end{cor}
\begin{proof} Under definition, we have $\theta \sb f \circ \alpha \sb {m(x)}=  \alpha \sb {f(x)}$. Apply $\Phi $
and obtain $\Phi (\theta \sb f )\circ \alpha \sb {m(x)}=  \Phi (\alpha \sb {f(x)})=\alpha \sb {s\sb A (f(x))}$.
Hence $\Q (\Phi (\theta \sb f))(m(x)) =s\sb A (f(x))$ for all $x\in X$.
\end{proof}

It turns out that the most interesting categories satisfy the condition we used above, namely, their
automorphisms take a free object over an one-element set to an isomorphic one. In such a case all automorphisms
are potential-inner and problem is to find a suitable extension for such a category. How to do it we show in the
next sections.

\section{ Generalized Reduction theorem}\label{sect_Reduct}
In this section we consider a category $\C$, its extension $\D$ with the same objects and a faithful functor $\Q
$ from $\D$ to the category of sets. Therefore $\D -$morphisms (and of course $\C -$morphisms) can be regarded
as a maps.

\begin{defn}\label{rightInd} An object $A\sp 0$ of the category $\C$ is called a {\it right indicator respectively }
$\D$ if for every two objects $A$ and $B$ and for every bijection $s: \Q (A )\to \Q (B )$ the following
condition is satisfied:

if for every $\C -$morphism $\nu: B \to A\sp 0  $  there exists a $\D -$morphism $\mu: A\to A \sp 0 $ such that
$ \Q (\mu )= \Q (\nu )\circ s $, then there exists a $\D -$isomorphism $\gamma : A \to B$ such that $\Q (\gamma
)=s$.
\end{defn}

Dually
\begin{defn}\label{leftInd}  An object $A\sp 0$ of the category $\C$ is called
a {\it left indicator respectively } $\D$ if for every two objects $A$ and $B$ and for every bijection $s: \Q (A
)\to \Q (B )$ the following condition is satisfied:

if for every $\C -$morphism $\nu: A\sp 0 \to A $ there exists a $\D -$morphism $\mu: A\sp 0 \to B$ such that $
\Q (\mu )= s\circ \Q (\nu ))$, then there exists a $\D -$isomorphism $\gamma : A \to B$ such that $\Q (\gamma
)=s$.
\end{defn}

Roughly speaking, the both conditions are the following ones: if composition of $s$ and a $\C -$morphism is a
$\D -$morphism then $s$ is a $\D -$isomorphism. In the case $\D = \C$ we use the term {\it indicator }. We give
some important examples of indicators in the varieties of universal algebras.

\begin{examples*} Let $\C$ be a category of universal algebras and all their homomorphisms. We choose in this case
the natural forgetful functor in the capacity of $\Q$.

\item[1.]  Let $A\sp 0$ be an algebra such that for every $\C$-object $A$ and every two its different elements
$a\sb 1$ and $a\sb 2$ there exists a homomorphism $\nu : A\to A\sp 0$ with $\nu (a\sb 1) \not = \nu (a\sb 2)$.
Then $A\sp 0$ is a right indicator in $\C$.

Indeed, let $s: A \to B$  be a bijection. Suppose that for every homomorphism $\nu: B \to A\sp 0 $
 the composition $\nu \circ s $ is a homomorphism too.  Consider an n-ary operation symbol $\omega$ and $n$ elements
$a\sb 1 ,\dots ,a\sb n \in A$. Let $a=\omega (a\sb 1 ,\dots ,a\sb n)$ and $b=\omega (s(a\sb 1) ,\dots ,s(a\sb
n))$. We have to show that $s(a)=b$. Suppose the contrary, that is,  $s(a)\not =b$. Under hypothesis there
exists a homomorphism $\nu : B \to A\sp 0$ with $\nu (s(a)) \not = \nu (b)$. We have: $(\nu \circ s)(a)=\omega
(\nu \circ s (a\sb 1),\dots ,\nu \circ s (a\sb n))=\omega (\nu (s(a\sb 1)),\dots ,\nu (s(a\sb n)))= \nu (\omega
(s(a\sb 1) ,\dots ,s(a\sb n)))=\nu (b).$

That means in contradiction to assumption that $\nu (s(a)) = \nu (b)$. Hence $s$ is a homomorphism and therefore
$A\sp 0$ is a right indicator.

\item[2.] Let $\C$ be a full subcategory of the category $\Theta \sp 0 (\V)$ for some variety $\V$ and let $F\sp
o$ be a free algebra generating the variety $\V$. It is easy to see that the algebra $A\sp 0 =F\sp 0$ satisfies
the conditions of the previous example. Thus if $F\sp 0$ is an object of $\C$ it is a right indicator in this
category.

\item[3.] Let $A\sp 0$ be a such object of $\C$ that for every $\C$-object $A$ and every finite subset $X$ of
$A$ having as many elements as arity of an operation, there exists a homomorphism $\nu : A\sp 0\to A$ with
$X\subset \nu (A\sp 0)$ . Then $A\sp 0$ is a left indicator in $\C$.

Indeed, let $s: A \to B$  be a bijection. Suppose that for every homomorphism $\nu:  A\sp 0 \to A$
 the composition $s \circ \nu $ is a homomorphism too.  Consider an n-ary operation symbol $\omega$ and $n$ elements
$a\sb 1 ,\dots ,a\sb n \in A$. Let $a=\omega (a\sb 1 ,\dots ,a\sb n)$ and $b=\omega (s(a\sb 1) ,\dots ,s(a\sb
n))$. We have to show that $s(a)=b$. Under hypothesis there exists a homomorphism $\nu :  A\sp 0 \to A $ with
$a\sb 1 ,\dots ,a\sb n \in \nu (A\sp 0 )$. It means that there exist $n$ elements $w\sb 1 ,\dots ,w\sb n \in
A\sp 0 $ such that $a\sb i = \nu (w\sb i)$ for all $i= 1,\dots , n$.  We have:
 $s(a)= s(\omega (a\sb 1 ,\dots ,a\sb n))=s(\omega (\nu (w\sb 1) ,\dots ,\nu (w\sb n))) =s\circ \nu (\omega (w\sb 1 ,\dots
,w\sb n))= \omega (s\circ \nu (w\sb 1),\dots ,s\circ \nu (w\sb n))=\omega (s(a\sb 1),\dots ,s(a\sb n))=b.$  That
means $s: A\to B $ is a homomorphism. And we make conclusion that $A\sp 0$ is a left indicator in $\C$.

\item[4.] As a corollary from the previous example, we obtain that every free algebra $H$ in $\C$ such that a
set of free generators of $H$ has not less elements than arities of all operations is a left indicator.
Particularly, every free algebra with two generators is a left indicator in a category of algebras with binary,
unary and nullary operations only.
\end{examples*}

The following result shows that if an automorphism of a category is indeed inner it is possible to detect this
considering only two objects.

\begin{thm}\label{main}Let $\C$ be a category with a forgetful functor $\Q: \C \to Set$ such that:
\begin{enumerate}
\item[1)]  functor $\Q $ is represented by an object $A_0$ ;

\item[2)] there is a right (or left) indicator  $A^0$ in $\C$.
\end{enumerate}
If $\Phi : \C \to \C $ is an automorphism of the category $\C $ that does not change the objects $A\sb 0 $ and $
A\sp 0 $ and induces the identity map on $Hom (A\sb 0 , A\sp 0) $ then $\Phi $ is an inner automorphism.
\end{thm}
\begin{proof}

Let us fix the objects $A\sp 0$ and $A\sb 0$ existing under hypotheses. Let $\Phi$ be an isomorphism of the
category $\C$ satisfying required conditions. According to \lemref{mainlem}, we define by \ref{E:sDefinition}
the family of bijections $(s\sb A :\Q (A)\to \Q(\Phi (A)),\; \vert A\in Ob\C )$, such that \ref{E:general} is
satisfied. Under hypothesis,  $\Phi $ lives fixed all morphisms from $A\sb 0 $ to $A\sp 0$, hence we have $s \sb
{A \sp 0}=1\sb {\Q (A\sp 0)}$.

Particularly, \ref{E:general} gives for every morphism $\nu :A\to A\sp 0 $ that $\Q (\Phi (\nu ))=s \sb {A\sp 0}
\circ \Q (\nu )\circ s \sb A \sp {-1}= \Q (\nu ) \circ s \sb A \sp {-1}$. Under hypotheses for the category
$\C$, there exists an isomorphism $\st \sb A :A\to \Phi (A)$ such that $s \sb A =\Q (\st \sb A )$. The dual case
gives the same conclusion. And finally we obtain that for every $\nu :A\to B $: $\Phi (\nu )=\st \sb B \circ \nu
\circ \st \sb A \sp {-1}$, that ends the proof.
\end{proof}

Now we combine together  \lemref{subcategory}, \lemref{objects} and  \thmref{main} to get a general result which
reduces the problem if an automorphism of a category is inner to more simple question.

\begin{thm}\label{general} Let functor $\Q $ be represented by an object $A_0$ and   $A^0$ be a left
(o right) indicator in $\C$. Then an automorphism $\Phi $ of $ \C$ is inner if and only if there are two
isomorphisms $\sigma : A\sb 0 \to \Phi (A\sp 0),\;\tau : A\sp 0\to \Phi (A\sp 0)$ such that for every
$\C$-morphism $\nu : A\sb 0\to A\sp 0$ we have $\Phi (\nu )= \tau \circ \nu \circ \sigma \sp {-1}$.
\end{thm}

\begin{proof}
The necessity of the given condition is obvious. We show that it is sufficient. Let $\Phi :\C \to \C$ is an
automorphism such that there are two isomorphisms $\sigma: A\sb 0 \to \Phi (A\sb 0),\;\tau: A\sp 0\to \Phi (A\sp
0)$ satisfying the mentioned conditions. According to \lemref{objects}, we can assume that $\Phi $ preserves
both mentioned objects. If $A\sb 0=A\sp 0$ then $\sigma =\tau $. Thus if objects $A_0$ and $A^0$ are different
and if they coincide we can apply \lemref{subcategory} and obtain that $\Phi $ is a composition of two
automorphisms $\Phi = \Psi \circ \Gamma $, where $\Gamma $ preserves objects $A\sb 0$ and $A\sp 0$ and all
morphisms $\nu :A \sb 0 \to A\sp 0$, and the functor $\Psi$ is an inner automorphism of the category $\C$.
According to \thmref{main}, $\Gamma $ is an inner automorphism of the category $\C$.
\end{proof}

The similar arguments lead to a  sufficient conditions for an automorphism to be $\D -$inner for a given
extension $\D $ of the category $\C$.

\begin{thm}\label{generalD-inner} Let functor $\Q $ be represented by an object $A_0$ and   $A^0$ be a left
(o right) indicator in $\C$ respectively $\D$. Let an automorphism $\Phi $ of $\C$ leave fixed these two
objects. Then  $\Phi$ is $\D -$inner if if the bijection $s\sb {A\sp 0}\sp {\Phi }$ is the $\Q -$value of a
$\D-$isomorphism.
\end{thm}

\begin{proof} Let  $A^0$ be a left indicator in $\C$ respectively $\D$. Since
$\Q (\Phi (\nu ))=s \sb A \circ \Q (\nu )\circ s \sb {A \sp 0} \sp {-1}$ for every $\C -$morphism $\nu : A\sp 0
\to A$, we have $s\sb A \circ \Q (\nu )=\Q (\Phi (\nu ))\circ s \sb {A \sp 0}$ for every $\C -$morphism $\nu :
A\sp 0 \to A$. We see that the right side of this equation is the $\Q -$value of a $\D -$morphism, because under
hypotheses $s\sb {A\sp 0}$ is the $\Q -$value od a $\D-$isomorphism. Under Definition~\ref{rightInd}, $s\sb A $
is the $\Q -$value of a $\D -$isomorphism $\sigma \sb A$ for every object $A$. Thus $\Q (\Phi (\nu ))=\Q (\sigma
\sb B \circ \nu \circ \sigma \sb A \sp {-1})$ for every $\C -$morphism $\nu : A\to B$. Since $\Q $ is faithful,
$\Phi $ is $\D$-inner. The same conclusion is clearly true in the case $A^0$ is a right indicator in $\C$
respectively $\D$.
\end{proof}

\section{ Categories of universal algebras}\label{sect_Univ}

In this section, we consider categories of universal algebras only. Given a variety $\V$ of universal algebras
of a type $\Xi$, $\Th (\V)$ denotes a the category of all $\V -$algebras and their homomorphisms. We apply our
results to a full subcategory $\C$  of the category $\Th (\V)$. In this case, a functor $\Q $ is the natural
forgetful functor. Therefore we denote homomorphisms and corresponding maps with the same letter, that is, we
write $\nu $ instead of $\Q (\nu ) $ for every homomorphism $\nu : A\to B $. To apply \thmref{general} it is
necessary to find two objects $A\sb 0 $ and $A\sp 0 $ satisfying conditions (1) and (2) respectively.  If a
monogenic free algebra $A\sb 0$ in $\Th (\V)$ is an object of $\C$ the condition (1) is realized. The condition
(2) is realized if our category contains an algebra $A\sp 0$ satisfying at least one of the conditions given in
examples in the previous section.

First we show how to obtain the original Reduction Theorem \cite{LieAlg}. Consider the category  $\C = \Th \sp
0$ described in Introduction, i.e.  $\Th \sp 0$ is the full subcategory of $\Th (\V )$ defined by all free
algebras from $\V$ over finite subsets of an infinite fixed set $X\sb 0$.  Let $F\sb 0$ be a free algebra in
$\Th \sp 0$ over an one-element set $\{x\sb 0\},\; x\sb 0 \in X\sb 0$.
\begin{thm}\label{Plotkin} \cite{LieAlg}
Assume that

\noindent (R1) every object of $\Th \sp 0 $ is a hopfian algebra and

\noindent (R2) there exists an object $F\sp 0=F(X)$ in $\Th \sp 0 $ generating the whole variety $\V$.

\noindent Let  $\nu \sb 0 :F\sp 0 \to F\sb 0 $ be the homomorphism induced by the constant map $X \to \{x\sb
0\}$, i.e. $ \nu \sb 0 (x)=x\sb 0$ for all $x\in X $. Under these conditions, if $\Phi : \Th \sp 0 \to \Th \sp
0$ is an automorphism such that 1) it does not change objects, 2) it induces the identity automorphism on the
semigroup $END (F\sp 0)$ and 3) it preserves $\nu \sb 0 :\; \Phi (\nu \sb 0 )=\nu \sb 0$, then $\Phi $ is an
inner automorphism.

\end{thm}
\begin{proof} We apply \thmref{main}. In the present case, $\C =\Th \sp 0 $ and $\Q :\Th \sp 0 \to Set $ is the
natural forgetful functor. Then the functor $\Q$ is represented by the object $F\sb 0$ and the condition (1) in
\thmref{main} for the category $\Th \sp 0 $ is satisfied. According to the Example~2 in Section~2, the algebra
$F\sp 0=F(X\sp 0 )$ is a right indicator in $\Th \sp 0 $, and therefore the condition~(2) is also satisfied.

Consider an automorphism $\Phi : \Th \sp 0 \to \Th \sp 0$.  Let $\mu : F\sb 0 \to F\sp 0$ be an arbitrary
homomorphism and let $\nu = \mu \circ \nu \sb 0 : F\sp 0 \to F \sp 0$. Under hypotheses, we have:
$$ \nu = \Phi (\nu )= \Phi (\mu )\circ \Phi (\nu \sb 0 )=\Phi (\mu ) \circ \nu \sb 0.$$
Hence $\mu \circ \nu \sb 0 = \Phi (\mu ) \circ \nu \sb 0$ and therefore $\Phi (\mu )= \mu $. Thus the
automorphism $\Phi$  satisfies the hypothesis of \thmref{main},  and therefore it is an inner automorphism.

\end{proof}

\begin{rem} The condition (R1) has not been used in the proof given above. Thus it is not necessary for the
\thmref{Plotkin}. However this condition is used in \cite{LieAlg} to prove the fact that every automorphism of
the category $\Th \sp 0 $ takes every object to an isomorphic one. But we have seen that it is also not
necessary, it is sufficient to assume this condition only for monogenic free algebra.

\end{rem}

Now we begin to characterize automorphisms which are not inner.

In most cases, every automorphism of the category $\Th (\V)$ takes the monogenic free algebras to isomorphic
algebras. This is an argument to assume further that an automorphism $\Phi $ of $\C$ preserves $A\sb 0$.
According to \ref{E:general}, for $a\in A $: $\Phi (\alpha \sb a)=s\sb A \circ \alpha \sb a \circ s\sb {A \sb
0}\sp {-1}\;$, where $s\sb {A\sb 0}$ is a permutation of $A\sb 0$.  Under definitions of maps $s\sb A$ above we
have $s\sb {A\sb 0}(x\sb 0)=x\sb 0$ and consequently for every $\alpha :A\sb 0 \to A$:
\begin{equation}\label{E:maps}
\Phi (\alpha )(x\sb 0)=(s\sb A \circ \alpha ) (x\sb 0).
\end{equation}

The bijections $s\sb A :A\to \Phi (A)$ can be used to define a new algebraic structure on the underlying set of
the algebra $ \Phi (A)$. We denote this new algebra by $A\sp *$, hence $s\sb A :A \to A\sp *$ is an isomorphism.
Of course, $A\sp * $ need not be an object of $\C$.

\begin{lem}\label{InnerAndCentral}
An automorphism $\Phi$ of $\C$ is inner if and only if there exists a central function $A\mapsto c\sb A$, $A\in
Ob\, \C$ such that  $c\sb {\Phi (A)}$ is an isomorphism of $\Phi (A)$ onto $A\sp *$.

\end{lem}
\begin{proof}
If $\Phi$ is inner, then there exists a function $A\mapsto \tau \sb A$, $A\in Ob\, \C$ such that  $\tau \sb A$
is an isomorphism of $A$ onto $\Phi (A)$ and for all $\nu :A \to B $ we have $\Phi (\nu )=\tau \sb B \circ \nu
\circ \tau \sb A \sp {-1 }$. Consider $c\sb {\Phi (A)}=s\sb A \circ \tau \sb A \sp {-1}$. Clearly, $c\sb {\Phi
(A)}$ is an isomorphism of $\Phi (A)$ onto $A\sp *$. On the other hand, for every homomorphism $\nu : \Phi
(A)\to \Phi (B)$ we have $c\sb {\Phi (B)} \circ \nu \circ c\sb {\Phi (A)} \sp {-1 }=s\sb B \circ (\tau \sb B \sp
{-1} \circ \nu \circ \tau \sb A )\circ s\sb A \sp {-1} = s\sb B\circ \Phi \sp {-1} (\nu )\circ s \sb A \sp {-1}
= \nu$. Hence the function $A\mapsto c\sb A$, $A\in Ob\, \C$ is central.

Inversely, if there exists such central function $A\mapsto c\sb A$, $A\in Ob\, \C$, that $c\sb {\Phi (A)}$ is an
isomorphism of $\Phi (A)$ onto $A\sp *$, then we set $\tau \sb A = c\sb {\Phi (A)} \sp {-1} \circ s\sb A $.
Clearly that $\tau \sb A$ is an isomorphism of $A$ onto $\Phi (A)$ and $\tau \sb B \circ \nu \circ \tau \sb A
\sp {-1 }= c\sb {\Phi (B)} \sp {-1}\circ s\sb B \circ \nu \circ s\sb A \sp {-1}\circ c\sb {\Phi (A)} =c\sb {\Phi
(B)} \sp {-1}\circ \Phi (\nu )\circ c\sb {\Phi (A)} =\Phi (\nu )$ for all $\nu :A \to B $. Hence $\Phi $ is
inner.
\end{proof}

The crucial fact is that in the case $A$ is a free algebra with enough number of free generators the structure
$A\sp *$ can be found and therefore the map $s\sb A$ can be described.

Denote by $r\sb \Xi $ the maximal arity of operations in $\Xi$. Let now $A$ be a free algebra which set of free
generators is $X=\{ x\sb 1 ,\ldots ,x\sb n \}$ where $n \geq r\sb \Xi$. Every element of $A$ being a term in
corresponding language determines  derivative (or polynomial, in other words) operations, which arities depend
on definition but are not greater than $n$. Particularly, suppose that $\Phi (A)=A$ and $\Phi (\alpha \sb x )=
\alpha \sb x $ for all $x\in X$. Let $\omega$ be a $k$-ary signature operation in $A$, $k\leq n$. Let $u= \omega
(x\sb 1,\ldots ,x\sb k )$ and $v$ be an element of $A$ defined by the equation: $\Phi (\alpha \sb u )=\alpha \sb
v $, that is, $v=s\sb A \sp {\Phi}(u)$. Then according to mentioned above $v$ determines a polynomial $k$-ary
operation $\omega \sp {\Phi } $ which can be expressed in the following way:
\begin{defn}\label{operation}
\begin{equation}\label{termoper}
\omega \sp {\Phi} (x\sb1 ,\ldots ,x\sb k )= s\sb A \sp {\Phi}(\omega (x\sb 1,\ldots ,x\sb k ))
\end{equation}
and for every elements $a\sb1 ,\ldots ,a\sb k $

\begin{equation}
\omega\sp {\Phi } (a\sb 1 ,\ldots ,a\sb k )=\theta \sb f ( \omega \sp {\Phi } (x\sb 1 ,\ldots ,x\sb k )) ,
\end{equation}
where $f(x\sb 1)=a\sb 1 ,\ldots ,f(x\sb k )=a\sb k ,f(x\sb {k+1})= x\sb {k+1}, \ldots ,f(x\sb n )=x\sb n $.
\end{defn}
Thus the automorphism $\Phi $  determines a new algebra $A\sp {\Phi } $  of the same type and with the same
underlying set that the algebra $A$ . We call this algebra a $\Phi$-derivative of the algebra $A$. It is clear
that along with $\Phi$-derivative we have $\Phi \sp {-1}$-derivative of the algebra.

\begin{thm}\label{main2}
Let $\C $ be a full subcategory of the category $\Th (\V)$ for some variety $\V$. Let $\C$ contain a monogenic
free algebra $A\sb 0$ and $A$ be a free finitely generated algebra in $\C$ with the set of free generators $X=\{
x\sb 1 ,\ldots ,x\sb n \} $ where $n \geq r\sb \Xi$. If $\Phi $ is an automorphism of  $\C $  that preserves
$A\sb 0$ and $A$, and $\Phi (\alpha \sb x )= \alpha \sb x $ for all $x\in X$, then $A\sp * =A\sp {\Phi }$.

\end{thm}
\begin{proof}
Denote $s=s\sb A $. Let $\omega $ be a $k$-ary signature operation of the algebra $A$ and $u= \omega (x\sb
1,\ldots ,x\sb k )$. Let $a\sb 1 ,\ldots ,a\sb k $ be elements of this algebra and $f(x\sb 1)=a\sb 1 ,\ldots
,f(x\sb k )=a\sb k ,f(x\sb {k+1})= x\sb {k+1}, \ldots ,f(x\sb n )=x\sb n $. Then  $\alpha \sb {\omega (s\sp {-1}
(a\sb 1 ),\ldots ,s\sp {-1}(a\sb k ))}=\theta \sb {s \sp {-1} \circ f } \circ \alpha \sb u$. Applying  $\Phi $
to this equation and using the definition $\omega \sp * (a\sb 1 ,\ldots ,a\sb k )= s(\omega (s \sp {-1} (a\sb 1
),\ldots ,s \sp {-1} (a\sb k )))$, we obtain:
$$ \alpha \sb {\omega \sp * (a\sb 1 ,\ldots ,a\sb k )}=\theta \sb  f \circ \alpha \sb {s(u)}
$$
or
$$ \omega \sp * (a\sb 1 ,\ldots ,a\sb k )=\omega \sp {\Phi} (a\sb 1 ,\ldots ,a\sb k ).
$$

\end{proof}
This result gives us an opportunity to reduce describing automorphisms of subcategories of $\Th (\V) $ to
studying derivative operations (polynomial operations) on free algebras of the variety $\V$. We have seen that
$A\sp *$ is derivative algebra of the same type than $A$ and is isomorphic to $A$, hence $A$ is derivative
algebra with respect $A\sp *$. It means that for every basic $k$-ary operation $\omega $, the value $\omega \sp
* (x\sb 1,\ldots ,x\sb k)$ is a polynomial $w$ in $A$ and $\omega (x\sb 1,\ldots ,x\sb k )$ is a polynomial $w
\sp *$ in $A\sp *$. Replacing in $w \sp *$ all operations of $A\sp *$ by its expressions as polynomials in $A$,
we obtain an identity $\omega (x\sb 1,\ldots ,x\sb k )=\overline{w \sp *} $ which must be satisfied in our
variety.

We show how these reasons simplify the problem, namely, we show how simple is to obtain some generalizations of
known results \cite{Free, InvSem} using our general approach.

\begin{thm}\label{semigroup}
Let $\C$ be a full subcategory of of the category of all semigroups containing a free monogenic semigroup $A\sb
0 = W(x\sb 0)$ and  a free semigroup $A\sp 0 = W(x,y)$ with two free generators $x,y$. Let the category $\D$ be
the extension of $\C$ by adding anti-homomorphisms. If an automorphism $\Phi$ of $\C$  takes $A\sb 0$
 and $A\sp 0 $ to isomorphic to them semigroups, then it is $\D$-inner.
\end{thm}
\begin{proof}
According to \lemref{objects}, we can assume that $\Phi$ preserves $A\sb 0$ and $A\sp 0$. Let $s=s\sb {A\sp 0
}$. According to \thmref{main2} we look for a derivative binary operation and corresponding semigroup $(A\sp 0
)\sp *$ with the same underlying set such that $s$ is an isomorphism preserving $x$ and $y$. Since the unique
identity of the kind $xy =w$ in variety of all semigroups is $xy=xy$, we conclude that there exist only two
derivative operations: $x\bullet y =xy$ and $x\ast y=yx$. Therefore $s(xy)=xy$ or $s(xy)=yx$. In the first case,
$s$ is the identity mapping. In the second one, $s$ maps every word $u$ to the "indirect" word $\bar{u}$, that
is all letters are written in reverse order. According to \thmref{generalD-inner}, $\Phi $ is $\D$-inner.
\end{proof}

It is proved in  \cite{InvSem} that the category of all free inverse semigroups is perfect. The proof uses a
description of $AUT END (F)$ for free inverse semigroups $F$. The next theorem generalizes this result using
\thmref{main2}.

We consider an inverse semigroup $A$ as an algebra with two operations, a binary multiplication $\cdot$ and a
unary inversion $\sp {-1}$ (here $a\sp {-1}$ is the inverse of an element $a$). The class of all inverse
semigroups forms a variety $\V $ defined by the identities:
$$(xy)z = x(yz),\;(xy)\sp {-1} = y\sp {-1} x\sp {-1},\;(x\sp {-1})\sp {-1} = x,$$
$$ xx\sp {-1}x = x,\; x\sp {-1} xy\sp {-1} y = y\sp {-1} yx\sp {-1} x.$$

\begin{thm}\label{inverse}
Let $\C$ be a full subcategory of of the category of inverse semigroups containing a free monogenic inverse
semigroup $A\sb 0 = W(x\sb 0)$ and  a free inverse semigroup $A\sp 0 = W(x,y)$ with two free generators $x,y$.
If an automorphism $\Phi$ of $\C$  takes $A\sb 0$  and $A\sp 0 $ to isomorphic to them semigroups, then it is
inner.
\end{thm}
\begin{proof}
The first step is the same that in \thmref{semigroup}, that is, we have a permutation $s$ of $A\sp 0 $  and we
look for a derivative inverse semigroup $(A\sp 0 )\sp *$ with the same underlying set such that $s$ is an
isomorphism of $A\sp 0 $ onto $(A\sp 0 )\sp *$, preserving $x$ and $y$. First of all, if a term $u$ determines
an involution it does not contain $y$ and thus $u=u (x)$. Hence we have the identity $x=u(u(x))$. It is known
that such identity is fulfilled in the variety of all inverse semigroups if and only if the word on the right
has the form $(xx\sp {-1})\sp k x$ or the form $x(x\sp {-1}x)\sp k$. Every part of these words is equal to $x$
or to $xx\sp {-1}$ or to $x\sp {-1}x$ or to $x\sp {-1}$. The first three words can not determine involutions.
Thus $u\equiv x\sp {-1}$ and $s(x\sp {-1}) =x\sp {-1}$.

Finally\footnote {The idea of the next step comes from G. Mashevitzky.}, let a term $w(x,y)$ determine a binary
operation and consider the following system of three equations: $w(x,x\sp {-1})=xx\sp {-1}$, $w (x,w (x\sp
{-1}x))=x$, $w(w(xx\sp {-1}),x)=x$. The only two terms satisfying this system are  $w=xy$ and $w=yx$. Thus
$s(xy)=xy$ or $s(xy)=yx$. In the first case, $s$ is the identity mapping and $(A\sp 0) \sp *$ coincides with
$A\sp 0$. In the second one, $(A\sp 0) \sp *$ is the dual inverse semigroup to $A\sp 0$. Because the involution
$c\sb A : a\mapsto a\sp {-1}$ is an isomorphism is an isomorphism of every inverse semigroup $A$ onto is dual
inverse semigroup $A\sp *$ and the function $A \mapsto c\sb A $ is a central, we conclide, according to
\lemref{InnerAndCentral}, that $\Phi $ is an inner automorphism.

\end{proof}

The similar method can be applied to categories of groups, modules, linear algebras and so on. In our next paper
\cite{Pl_Zh}, we apply this method and characterize automorphisms of categories  $\Th \sp 0 (\V) $ where $\V $
is the variety of all associative $K-$algebras and where $\V $ is the variety of all free group representations.

\section {Some other applications}\label{applications}

We should mention that potential-inner automorphisms in fact are inner in the case of full subcategories of the
category $Set$. Thus the following result is a trivial corollary of \thmref{potent-inner}.

\begin{thm}\label{set} Every full subcategory of the category  $Set$ containing one-element set is perfect.
\end{thm}
The following old classical result immediately follows from \thmref{potent-inner}, in other words,
\thmref{potent-inner} can be regarded as a generalization of this result to categories.
\begin{thm}\label{Schreier}\cite{Schreier} Every automorphism of the semigroup $T \sb X $ of all transformations of some
set $X$ is inner.
\end{thm}
\begin{proof}Consider the category with one object $X$ and the set $S \sb X $ as the set of all its morphisms.
Let functor $\Q$ be $Hom (X, -)$. It is faithful and representable under definition. According to
\lemref{mainlem}, all automorphisms of this category are inner.
\end{proof}

The next also known old result follows from \thmref{set}.
\begin{thm}\label{Shutov}\cite{Shutov} Every automorphism of the semigroup $\mathfrak F \sb X $ of all partial
transformations of a set $X$ is inner.
\end{thm}
\begin{proof}
Consider the category $\mathcal P (X)$ of all subsets of $X$ and their mappings. According to \thmref{set} it is
perfect. Denote by $\Delta \sb A $ the identity relation on the set $ A\subseteq X$: $\Delta \sb A =\{(a,a)\vert
a\in A\}$. Let $\Phi $ be an automorphism of the semigroup $\mathfrak F \sb X $. It is known that $\Phi (\Delta
\sb A )=\Delta \sb B$ for some $ B\subseteq X$. Thus we define $\Phi ( A )= B \Leftrightarrow \Phi (\Delta \sb A
)=\Delta \sb B$ and consider $\Phi$ as an automorphism of the category $\mathcal P (X)$. Since it is inner there
is a family of bijections $(s\sb A :A\to \Phi (A) \vert  A\subseteq X )$  such that for every partial
transformation $f$ considered as a map $f:A\to B$ the following equality takes place: $\Phi (f)=s\sb B \circ f
\circ s\sb A \sp {-1}$. Particularly for $f=\Delta \sb A :A\to X$, we have $\Delta \sb {\Phi (A)}=\Phi (\Delta
\sb A )=s\sb X \circ \Delta \sb A \circ s\sb A \sp {-1}=s\sb X \circ  s\sb A \sp {-1}$. The last means that $s
\sb A $ is the restriction of $s\sb X $ on $A$. Thus $\Phi (f)=s\sb X \circ f \circ s\sb X \sp {-1}$ for every
$f\in \mathfrak F \sb X $.

\end{proof}
The same reasons are valid for the inverse semigroup of all one-to-one partial transformations of a given set
$X$.

Institute of Mathematics, Hebrew University, 91803 Jerusalem, Israel

{\it E-mail address}: borisov@math.huji.ac.il

Department of Mathematics, Bar-Ilan University, 52900 Ramat Gan, Israel

{\it E-mail address}: zhitomg@macs.biu.ac.il

\end{document}